# HARNACK INEQUALITY AND APPLICATIONS FOR STOCHASTIC GENERALIZED POROUS MEDIA EQUATIONS[1]

By Feng-Yu Wang

*Beijing Normal University*

By using coupling and Girsanov transformations, the dimension-free Harnack inequality and the strong Feller property are proved for transition semigroups of solutions to a class of stochastic generalized porous media equations. As applications, explicit upper bounds of the $L^p$-norm of the density as well as hypercontractivity, ultracontractivity and compactness of the corresponding semigroup are derived.

**1. Introduction.** The dimension-free Harnack inequality, first introduced by the author in [19] for diffusions on Riemannian manifolds, has been applied and extended intensively in the study of finite- and infinite-dimensional diffusion semigroups; see, for example, [16, 17, 20, 22] for applications to contractivity properties and functional inequalities, [1, 2, 11] for applications to short-time behaviors of infinite-dimensional diffusions, and [7, 8] for applications to the transportation-cost inequality and heat kernel estimates.

To establish the dimension-free Harnack inequality, the gradient estimate of the type $|\nabla P_t f| \leq e^{Kt} P_t |\nabla f|$ has played a key role in the above mentioned references, where the gradient is induced by the underlying diffusion coefficient. On the other hand, however, in many cases the semigroup is not regular enough to satisfy this gradient estimate; indeed, this gradient estimate is equivalent to Bakry–Emery's curvature condition for a very general framework as in [5]. To establish the dimension-free Harnack inequality on manifolds with unbounded below curvatures, a new approach is developed in the recent work [3] by using coupling and Girsanov transformations.

Received October 2005; revised June 2006.

[1]Supported in part by National Natural Science Foundation of China Grant 10121101 and Rural Finance Development Program Grant 20040027009.

In this paper, we intend to study the transition semigroup for solutions to a class of stochastic generalized porous media equations, for which the semigroup is merely known to be Lipschitzian in the natural norm rather than in the intrinsic distance (cf. [6]). So, we are not able to prove the Harnack inequality by using intrinsic gradient estimates. On the other hand, since the intrinsic distance is usually too big to be exponential integrable w.r.t. the underlying reference measure, we prefer to establish a Harnack inequality depending only on the natural norm. Such a stronger inequality will provide more information including the strong Feller property and the ultracontractivity of the semigroup. To modify the argument in [3], we shall construct a new coupling which only depends on the natural distance rather than the intrinsic one between the marginal processes (see Section 2 below).

Strong solutions of the stochastic generalized porous medium equation have been studied intensively in recent years; see [6] for the existence, uniqueness and long-time behavior of some stochastic generalized porous media equations with finite reference measures, see [12] for the stochastic porous media equation on $\mathbf{R}^d$ where the reference (Lebesgue) measure is infinite; and see [18] for large deviation principles. Recently, a general result concerning existence and uniqueness was presented in [15] for strong solutions of stochastic generalized porous media and fast diffusion equations.

Let $(E, \mathcal{M}, \mathbf{m})$ be a separable probability space and $(L, \mathcal{D}(L))$ a negative definite self-adjoint linear operator on $L^2(\mathbf{m})$ having discrete spectrum. Let

$$(0 <)\lambda_1 \leq \lambda_2 \leq \cdots$$

be all eigenvalues of $-L$ with unit eigenfunctions $\{e_i\}_{i \geq 1}$.

To state our equation, we first introduce the state space of the solutions. Let $H$ be the completion of $L^2(\mathbf{m})$ under the inner product

$$\langle x, y \rangle_H := \sum_{i=1}^{\infty} \frac{1}{\lambda_i} \langle x, e_i \rangle \langle y, e_i \rangle,$$

where $\langle \cdot, \cdot \rangle$ is the inner product in $L^2(\mathbf{m})$. It is well known that $H$ is the dual space of the Sobolev space $H^1 := \mathcal{D}((-L)^{1/2})$ and hence, is often denoted by $H^{-1}$ in the literature. Let $\mathcal{L}_{\mathrm{HS}}$ denote the space of all Hilbert–Schmidt operators from $L^2(\mathbf{m})$ to $H$. Let $W_t$ be the cylindrical Brownian motion on $L^2(\mathbf{m})$ w.r.t. a complete filtered probability space $(\Omega, \mathcal{F}_t, \mathbf{P})$; that is, $W_t := \{B_t^i e_i\}_{i \geq 1}$ for a sequence of independent one-dimensional $\mathcal{F}_t$-Brownian motions $\{B_t^i\}$. Let

$$\Psi, \Phi : [0, \infty) \times \mathbf{R} \times \Omega \to \mathbf{R}$$

be progressively measurable and continuous in the second variable, and let

$$Q : [0, \infty) \times \Omega \to \mathcal{L}_{\mathrm{HS}}$$

be progressively measurable such that

$$\mathbf{E}\int_0^T \|Q_t\|_{\mathcal{L}_{\mathrm{HS}}}^2 \, dt < \infty, \quad T > 0. \tag{1.1}$$

We consider the equation

$$dX_t = \{L\Psi(t, X_t) + \Phi(t, X_t)\} \, dt + Q_t \, dW_t. \tag{1.2}$$

In particular, if $\Phi = 0, Q = 0$ and $\Psi(t, s) := |s|^{r-1}s$ for some $r > 1$, then (1.2) reduces back to the classical porous medium equation (see, e.g., [4]).

In general, for a fixed number $r \geq 1$, we assume that there exist functions $\delta, \eta, \gamma, \sigma \in C([0, \infty))$ with $\delta > 0$ such that

$$\begin{aligned}
&|\Psi(t,s)| + |\Phi(t,s) - \sigma_t s| \leq \eta_t(1 + |s|^r), \qquad s \in \mathbf{R}, \, t \geq 0, \\
&2\langle \Psi(t,x) - \Psi(t,y), y - x \rangle - 2\langle \Phi(t,x) - \Phi(t,y), L^{-1}(x-y)\rangle \\
&\qquad \leq -\delta_t^2 \|x-y\|_{r+1}^{r+1} + \gamma_t \|x-y\|_H^2, \qquad x, y \in L^{r+1}(\mathbf{m}), \, t \geq 0,
\end{aligned} \tag{1.3}$$

where and in the sequel, $\|\cdot\|_p$ denotes the norm in $L^p(\mathbf{m})$ for $p \geq 1$. A very simple example satisfying (1.3) is that $\Psi(t, s) := |s|^{r-1}s$ and $\Phi(t, s) := \gamma_t s$.

By the first inequality in (1.3), the first term in the left-hand side of the second inequality makes sense for any $x, y \in L^{r+1}(\mathbf{m})$. Since $L^{-1}$ is bounded in $L^2(\mathbf{m})$, if $|\Phi(t,s)| \leq \sigma_t(1 + |s|^{(r+1)/2})$ for some positive $\sigma \in C([0,\infty))$, then the another term $\langle \Phi(t,x) - \Phi(t,y), L^{-1}(x-y)\rangle$ makes sense too. Otherwise, since the first condition in (1.3) only implies $|\Phi(t,s)| \leq \eta_t(1 + |s|^r)$, in general, to make the second condition in (1.3) meaningful, we should and do assume that $L^{-1}$ is bounded in $L^{r+1}(\mathbf{m})$. In particular, this assumption holds automatically if $L$ is a Dirichlet operator (cf., e.g., [14]).

Recall that an adapted continuous process $X_t$ is called a solution to (1.2) if (cf. [6])

$$\mathbf{E}\int_0^T \|X_t\|_{r+1}^{r+1} \, dt < \infty, \qquad T > 0,$$

and for any $f \in L^{r+1}(\mathbf{m})$,

$$\langle X_t, f \rangle_H = \langle X_0, f \rangle_H - \int_0^t \mathbf{m}(f\Psi(s, X_s) + \Phi(s, X_s)L^{-1}f) \, ds$$
$$+ \int_0^t \langle Q(s, X_s) \, dW_s, f \rangle_H, \qquad t \geq 0.$$

Due to (1.1), (1.3) and Theorems II.2.1 and II.2.2 in [13], for any $X_0 \in L^2(\Omega \to H; \mathcal{F}_0, \mathbf{P})$ the equation (1.2) has a unique solution (cf. Theorem A.2 below). For any $x \in H$, let $X_t(x)$ be the unique solution to (1.2) with $X_0 = x$. Define

$$P_t F(x) := \mathbf{E} F(X_t(x)), \qquad x \in H,$$

for any bounded measurable function $F$ on $H$.

We first study Harnack inequalities for $P_t$. To this end, we assume that $Q_t(\omega)$ is nondegenerate for $t > 0$ and $\omega \in \Omega$; that is, $Q_t(\omega)x = 0$ implies $x = 0$. Let

$$\|x\|_{Q_t} := \begin{cases} \|y\|_2, & \text{if } y \in L^2(\mathbf{m}), Q_t y = x, \\ \infty, & \text{otherwise.} \end{cases}$$

We call $\|\cdot\|_{Q_t}$ the intrinsic distance induced by $Q_t$.

THEOREM 1.1. *Assume* (1.1) *and* (1.3). *If there exists a nonnegative constant* $\theta \in r - 3$ *such that*

$$(1.4) \qquad \|x\|_{r+1}^{r+1} \geq \xi_t^2 \|x\|_{Q_t}^{2+\theta} \|x\|_H^{r-1-\theta}, \qquad x \in L^{r+1}(\mathbf{m}),\ t \geq 0,$$

*holds on* $\Omega$ *for some strictly positive function* $\xi \in C([0, \infty))$, *then for any* $t > 0$, $P_t$ *is strong Feller and for any positive bounded measurable function* $F$ *on* $H$, *any* $\alpha > 1$ *and any* $x, y \in H$,

$$(1.5) \qquad (P_t F)^\alpha(y) \leq (P_t F^\alpha(x)) \exp\left[\frac{\alpha c(\theta, t) \|x - y\|_H^{2(3-r+\theta)/(2+\theta)}}{(\alpha - 1)}\right],$$

*where*

$$c(\theta, t) := \left(2(4+\theta)^{(6+2\theta)/(2+\theta)} \left(\int_0^t \delta_s^2 \xi_s^2 \exp\left[-\frac{3-r+\theta}{4+\theta}\int_0^s \gamma_u\, du\right] ds\right)^{\theta/(2+\theta)}\right)$$

$$\times \left((3-r+\theta)^{(6+2\theta)/(2+\theta)}\right.$$

$$\left.\times \left(\int_0^t \delta_s \xi_s \exp\left[-\frac{3-r+\theta}{4+\theta}\int_0^s \gamma_u\, du\right] ds\right)^2\right)^{-1}.$$

Unlike known Harnack inequalities established in [1, 2, 11] where the involved distance is almost surely infinite, (1.5) only includes the usual norm on the state space $H$. This enables one to derive stronger regularity properties of the semigroup, such as the strong Feller property of $P_t$ and estimates of its transition density $p_t(x, y)$. Moreover, as was done in [16, 19, 20], this inequality can also be applied to derive the hypercontractivity and ultracontractivity of the semigroup (cf. Theorem 1.2 below).

To apply Theorem 1.1 to contractivity properties of $P_t$, we consider the following time-homogenous case.

THEOREM 1.2. *Assume* (1.1), (1.3) *and* (1.4) *for some nonnegative constant* $\theta > r - 3$. *Furthermore, let* $\Psi, \Phi$ *and* $Q$ *be deterministic and time-free such that* $\xi, \delta > 0$ *and* $\gamma$ *are constant with* $\gamma 1_{\{r=1\}} < \delta^2 \lambda_1$.

(1) *The Markov semigroup $P_t$ has an invariant probability measure $\mu$ with full support on $H$ and $\mu(e^{\varepsilon_0 \|\cdot\|_H^{r+1}} + \|\cdot\|_{r+1}^{r+1}) < \infty$ for some $\varepsilon_0 > 0$. If in addition $\gamma \leq 0$, then the invariant probability measure is unique.*

(2) *For any $x \in H$, any $t > 0$ and any $\alpha > 1$, the transition density $p_t(x,y)$ of $P_t$ w.r.t. $\mu$ satisfies*

$$\begin{aligned}(1.6) \quad &\|p_t(x,\cdot)\|_{L^p(\mu)} \\ &\leq \Bigg\{\int_H \exp\bigg[-(\alpha c(\theta)\|x-y\|_H^{2(3-r+\theta)/(2+\theta)}) \\ &\qquad \times \bigg(\bigg\{\frac{1}{\gamma}\bigg(1 - \exp\bigg[-\frac{3-r+\theta}{4+\theta}\gamma t\bigg]\bigg)\bigg\}^{(4+\theta)/(2+\theta)}\bigg)^{-1}\bigg]\mu(dy)\Bigg\}^{-(\alpha-1)/\alpha},\end{aligned}$$

*where $c(\theta) := 2(4+\theta)/(3-r+\theta)(\xi\delta)^{4/(4+\theta)}$ and when $\gamma = 0$, the right-hand side means its limit as $\gamma \downarrow 0$.*

(3) *If $r = 1$, then $P_t$ is hyperbounded (i.e., $\|P_t\|_{L^2(\mu) \to L^4(\mu)} < \infty$) and compact on $L^2(\mu)$ for some $t > 0$. If moreover $\gamma \leq 0$, then $P_t$ is hypercontractive, that is, $\|P_t\|_{L^2(\mu) \to L^4(\mu)} \leq 1$ for large $t > 0$.*

(4) *If $r > 1$, then $P_t$ is ultracontractive and compact on $L^2(\mu)$ for any $t > 0$. More precisely, there exists $c > 0$ such that*

$$(1.7) \qquad \|P_t\|_{L^2(\mu) \to L^\infty(\mu)} \leq \exp[c(1 + t^{-(1+r)/(r-1)})], \qquad t > 0.$$

To apply Theorems 1.1 and 1.2, one has to verify condition (1.4). To this end, we present below some simple sufficient conditions for (1.4) to hold.

COROLLARY 1.3. *Let $Qe_i := q_i e_i$ for $i \geq 1$ with $\sum_{i=1}^\infty \frac{q_i^2}{\lambda_i} < \infty$, so that $Q$ is Hilbert–Schmidt from $L^2(\mathbf{m})$ to $H$. If $\inf_i q_i^2 > 0$, then (1.4) holds for any nonnegative constant $\theta \in (r-3, r-1]$ and a constant function $\xi > 0$. Consequently, if moreover $\Psi$ and $\Phi$ are deterministic and time-free such that (1.3) holds with $\gamma 1_{\{r=1\}} < \lambda_1 \delta^2$, then all assertions in Theorems 1.1 and 1.2 hold for $\theta \in (r-3, r-1] \cap [0, \infty)$.*

PROOF. Simply note that $\|\cdot\|_{r+1}^2 \geq \|\cdot\|_2^2 \geq \frac{1}{\inf_i q_i^2}\|\cdot\|_Q^2$. □

REMARK 1.1. In Corollary 1.3 there are two conditions on $q_i$, where $\sum_{i \geq 1} \frac{q_i^2}{\lambda_i} < \infty$ means that $\{q_i^2\}$ should be small enough as $i \to \infty$ but the other says that the sequence should be at least uniformly positive. In particular, such sequence exists if the spectrum of $L$ is discrete enough such that $\sum_{i \geq 1} \frac{1}{\lambda_i} < \infty$. This is the case if, for example, $L = \Delta$ on a bounded domain in $\mathbf{R}$ with the Dirichlet boundary condition, or more generally, $L$ is the Laplace operator on a post-critical finite self-similar fractal with $s > 0$ the Hausdorff dimension of the fractal in the effective resistance metric (see [10]). In the first case it is well known that $\lambda_i \geq ci^2$ for some $c > 0$ and all $i \geq 1$, while according to Theorem 2.11, for the second case one has $\lambda_i \geq ci^{(s+1)/s}$ for some $c > 0$ and all $i \geq 1$. See Section 3 below for more examples of $L$ in an abstract framework including high-order elliptic differential operators on $\mathbf{R}^d$.

Complete proofs of the above two theorems will be presented in Section 2. Assertions in Theorem 1.2 are direct consequences of Theorem 1.1 as soon as the desired concentration of $\mu$ is confirmed. To prove the first theorem, we adopt the coupling method and Girsanov transformations as in [3]. Comparing to the argument developed in [19], this method enables one to avoid verifying (intrinsic) gradient estimates of the semigroup.

In Section 3, concrete sufficient conditions for Corollary 1.3 to hold are provided for a large class of linear operators $L$ in a rather abstract framework. Finally, in the Appendix we confirm the existence and uniqueness of the solution to (1.2) as well as the existence and uniqueness of our coupling constructed below [cf. (2.2)].

## 2. Proofs of Theorems 1.1 and 1.2.

2.1. *The main idea.* To make the proofs easy to follow, let us first briefly explain the main idea to obtain a Harnack inequality using coupling. Let $x \neq y$ be two fixed points in $H$, and let $T > 0$ be a fixed time. Let $X_t(x)$ and $X_t(y)$ be the solutions to (1.2) with initial data $x$ and $y$, respectively. If

$$(2.1) \qquad \tau(x,y) := \inf\{t \geq 0 : X_t(x) = X_t(y)\} \leq T \qquad \text{a.s.},$$

then by the uniqueness of the solution, we have $X_T(x) = X_T(y)$ a.s. Thus, for any nonnegative measurable function $F$ on $H$,

$$P_T f(x) := \mathbf{E} F(X_T(x)) = \mathbf{E} F(X_T(y)) = P_T F(y).$$

This is much more than the Harnack inequality we wanted. Of course, in general (2.1) is wrong since it is so strong that $P_T$ maps any bounded function to constant. What we can hope is that $\tau(x,y) \leq T$ happens in a high probability (for $x$ and $y$ close enough). This is, however, not sufficient to imply the Harnack inequality.

To ensure that $\tau(x,y) \leq T$ happens in probability 1, we shall add a strong enough drift term which forces $X_t(y)$ to move to $X_t(x)$. To this end, let us take a constant $\varepsilon \in (0,1)$ and a reference function $\beta \in C([0,\infty); \mathbf{R}_+)$, and consider the modified equation

$$(2.2) \quad dY_t = \left\{ L\Psi(t, Y_t) + \Phi(t, Y_t) + \frac{\beta_t(X_t - Y_t)}{\|X_t - Y_t\|_H^\varepsilon} 1_{\{t < \tau\}} \right\} dt + Q_t \, dW_t,$$

$$Y_0 = y,$$

where $X_t := X_t(x)$ and $\tau := \inf\{t \geq 0 : X_t = Y_t\}$.

By Theorem A.2 below, (2.2) has a unique solution. Moreover, by the uniqueness, we have $X_t = Y_t$ for $t \geq \tau$.

Now, to derive the desired Harnack inequality, we need only to find out $\varepsilon > 0$ and nonnegative function $\beta_t$ such that:

  (i) $\tau \leq T$ a.s.
  (ii) $\mathbf{E} \exp[\int_0^T \frac{\beta_t^2}{2} \|X_t - Y_t\|_H^{-2\varepsilon} \|X_t - Y_t\|_{Q_t}^2 \, dt] < \infty$.

Let
$$\zeta_t := \frac{\beta_t Q_t^{-1}(X_t - Y_t)}{\|X_t - Y_t\|_H^\varepsilon} 1_{\{t < \tau\}}.$$

Once (i) and (ii) are confirmed, we may rewrite (2.2) as
$$dY_t = (L\Psi(t, Y_t) + \Phi(t, Y_t)) \, dt + Q_t \, d\tilde{W}_t, \qquad Y_0 = y,$$

where
$$\tilde{W}_t := W_t + \int_0^t \zeta_s \, ds, \qquad t \in [0, T].$$

By (ii) and Girsanov's theorem, it is easy to see that $\{\tilde{W}_t\}_{t \in [0,T]}$ is a cylindrical Brownian motion on $L^2(\mathbf{m})$ under the weighted probability measure $R\mathbf{P}$, where
$$R := \exp\left[\int_0^T \langle dW_t, \zeta_t \rangle - \tfrac{1}{2} \int_0^T \|\zeta_t\|_2^2 \, dt\right].$$

Thus, by the uniqueness of the solution, the distribution of $\{Y_t\}_{t \in [0,T]}$ under $R\mathbf{P}$ coincides with that of $\{X_t(y)\}_{t \in [0,T]}$ under $\mathbf{P}$. Therefore, combining this with (i) we arrive at

$$\begin{aligned}
P_T F(y) = \mathbf{E} R F(Y_T) &= \mathbf{E} R F(X_T) \\
&\leq (\mathbf{E} R^{\alpha/(\alpha-1)})^{(\alpha-1)/\alpha} (\mathbf{E} F(X_T)^\alpha)^{1/\alpha} \\
&= (\mathbf{E} R^{\alpha/(\alpha-1)})^{(\alpha-1)/\alpha} (P_T F^\alpha(x))^{1/\alpha}.
\end{aligned} \quad (2.3)$$

Then the desired Harnack inequality follows by estimating the moments of $R$.

2.2. *Proofs.* We first study (i). By (1.3) and the Itô formula due to [13], Theorem I.3.2, we have

$$\begin{aligned}
d\|X_t - Y_t\|_H^2 \\
\leq (-\delta_t^2 \|X_t - Y_t\|_{r+1}^{r+1} + \gamma_t \|X_t - Y_t\|_H^2 - \beta_t \|X_t - Y_t\|_H^{2-\varepsilon}) \, dt, \qquad t \leq T.
\end{aligned}$$

Then

$$d\{\|X_t - Y_t\|_H^2 e^{-\int_0^t \gamma_s ds}\}$$
(2.4)
$$\leq -(\delta_t^2 \|X_t - Y_t\|_{r+1}^{r+1} + \beta_t \|X_t - Y_t\|_H^{2-\varepsilon}) e^{-\int_0^t \gamma_s ds} dt, \qquad t \leq T.$$

LEMMA 2.1. *If $\beta$ satisfies*

$$\int_0^T \exp\left[-\frac{\varepsilon}{2} \int_0^t \gamma_s \, ds\right] \beta_t \, dt \geq \frac{2}{\varepsilon} \|x - y\|_H^\varepsilon, \tag{2.5}$$

*then $X_T = Y_T$.*

PROOF. By (2.4),

$$\frac{2}{\varepsilon} d\{\|X_t - Y_t\|_H^2 e^{-\int_0^t \gamma_s ds}\}^{\varepsilon/2} \leq -\beta_t e^{-\varepsilon/2 \int_0^t \gamma_s ds} dt, \qquad t \leq \tau.$$

If $T < \tau$, then it follows from this and (2.5) that

$$\{\|X_T - Y_T\|_H^2 e^{-\int_0^T \gamma_s ds}\}^{\varepsilon/2} - \|x - y\|_H^\varepsilon \leq -\frac{\varepsilon}{2} \int_0^T \beta_t e^{-\varepsilon/2 \int_0^t \gamma_s \, ds} \, dt$$
$$\leq -\|x - y\|_H^\varepsilon.$$

This implies $X_T = Y_T$ and hence, is contradictory to $T < \tau$. □

PROOF OF THEOREM 1.1. By (2.4), (1.4) and letting $\varepsilon := (3-r+\theta)/(4+\theta)$ which is in $(0,1)$ since $\theta > r - 3$, we obtain

(2.6)
$$d\{\|X_t - Y_t\|_H^2 e^{-\int_0^t \gamma_s ds}\}^\varepsilon$$
$$\leq -\varepsilon \delta_t^2 \|X_t - Y_t\|_H^{2(\varepsilon-1)} e^{-\varepsilon \int_0^t \gamma_s ds} \|X_t - Y_t\|_{r+1}^{r+1} dt$$
$$\leq -\varepsilon \delta_t^2 \xi_t^2 \|X_t - Y_t\|_{Q_t}^{2+\theta} e^{-\varepsilon \int_0^t \gamma_s ds} \|X_t - Y_t\|_H^{2(\varepsilon-1)+r-1-\theta} dt$$
$$= -\varepsilon \delta_t^2 \xi_t^2 e^{-\varepsilon \int_0^t \gamma_s ds} \frac{\|X_t - Y_t\|_{Q_t}^{2+\theta}}{\|X_t - Y_t\|_H^{(2+\theta)\varepsilon}} dt.$$

Let

$$\beta_t^2 := c^2 \delta_t^2 \xi_t^2 e^{-\varepsilon \int_0^t \gamma_s ds}, \qquad c := \frac{2\|x - y\|_H^\varepsilon}{\varepsilon \int_0^T \delta_t \xi_t \exp[-\varepsilon \int_0^t \gamma_s] \, ds}. \tag{2.7}$$

Then (2.5) holds so that $X_T = Y_T$ according to Lemma 2.1. So, (2.6) implies

$$\frac{\varepsilon}{c^2} \int_0^T \frac{\beta_t^2 \|X_t - Y_t\|_{Q_t}^{2+\theta}}{\|X_t - Y_t\|_H^{(2+\theta)\varepsilon}} dt \leq \|x - y\|_H^{2\varepsilon}.$$

STOCHASTIC POROUS MEDIUM EQUATION 9By this and the Hölder inequality,

$$
(2.8) \quad \begin{aligned}
\int_0^T & \frac{\beta_t^2 \|X_t - Y_t\|_{Q_t}^2}{\|X_t - Y_t\|_H^{2\varepsilon}} \, dt \\
& \leq \left( \int_0^T \frac{\beta_t^2 \|X_t - Y_t\|_{Q_t}^{2+\theta}}{\|X_t - Y_t\|_H^{(2+\theta)\varepsilon}} \, dt \right)^{2/(2+\theta)} \left( \int_0^T \beta_t^2 \, dt \right)^{\theta/(2+\theta)} \\
& \leq (\varepsilon^{-1} c^2 \|x - y\|_H^{2\varepsilon})^{2/(2+\theta)} \left( \int_0^T \beta_t^2 \, dt \right)^{\theta/(2+\theta)}.
\end{aligned}
$$

This implies, for $\alpha' := \alpha/(\alpha - 1)$, that

$$
(2.9) \quad \begin{aligned}
\mathbf{E} R^{\alpha'} &= \mathbf{E} \exp\left[ \alpha' \int_0^T \langle dW_t, \zeta_t \rangle - \frac{\alpha'}{2} \int_0^T \|\zeta_t\|_2^2 \, dt \right] \\
&= \mathbf{E} \exp\left[ \frac{\alpha'(\alpha' - 1)}{2} \int_0^T \|\zeta_t\|_2^2 \, dt \right] \\
&\leq \exp\left[ \frac{\alpha'(\alpha' - 1)}{2} (\varepsilon^{-1} c^2 \|x - y\|_H^{2\varepsilon})^{2/(2+\theta)} \left( \int_0^T \beta_t^2 \, dt \right)^{\theta/(2+\theta)} \right].
\end{aligned}
$$

Combining (2.9) with (2.3), we arrive at

$$
(P_T F(y))^\alpha \\
\leq (P_T F^\alpha)(x) \exp\left[ \frac{\alpha}{2(\alpha - 1)} (\varepsilon^{-1} c^2 \|x - y\|_H^{2\varepsilon})^{2/(2+\theta)} \left( \int_0^T \beta_t^2 \, dt \right)^{\theta/(2+\theta)} \right].
$$

Taking (2.7) into account, we obtain (1.5).

We now prove the strong Feller property. Since

$$
P_T F(y) = \mathbf{E} R F(Y_T) = \mathbf{E} R F(X_T),
$$

we have

$$
(2.10) \quad |P_T F(y) - P_T F(x)| = |\mathbf{E}(R - 1) F(X_T)| \leq \|F\|_\infty \mathbf{E} |R - 1|.
$$

From (2.9) we know that $R$ is uniformly integrable for bounded $\|x - y\|_H$. Therefore, by (2.8) and the dominated convergence theorem we obtain

$$
\lim_{y \to x} \mathbf{E} |R - 1| = \mathbf{E} \lim_{y \to x} |R - 1| = 0.
$$

Combining this with (2.10) we see that $P_T F \in C_b(H)$. Thus, $P_T$ is strong Feller. □

PROOF OF THEOREM 1.2(1). (a) The existence of $\mu$. Let $X_t(0)$ be the solution to (1.2) with $X_0 = 0$, and let

$$
\mu_n := \frac{1}{n} \int_0^n \delta_0 P_t \, dt, \qquad n \geq 1,
$$

where $\delta_0 P_t$ is the distribution of $X_t(0)$, $t \geq 0$. Since by Theorem 1.1 $P_t$ is a (even strong) Feller Markov semigroup, to prove the existence of the invariant probability measure, we only need to verify the tightness of $\{\mu_n : n \geq 1\}$. Indeed, if $\mu_{n_k} \to \mu$ weakly for some subsequence $n_k \to \infty$, then for any $F \in C_b(H)$ one has $P_t F \in C_b(H)$ and thus,

$$(\mu P_t)(F) = \lim_{k \to \infty} \mu_{n_k}(P_t F) = \lim_{k \to \infty} \frac{1}{n_k} \int_t^{n_k+t} P_s F(0)\, ds$$
$$= \lim_{k \to \infty} \frac{1}{n_k} \int_0^{n_k} P_s F(0)\, ds = \mu(F), \qquad t \geq 0.$$

By (1.3) with $\delta > 0$ and $\gamma 1_{\{r=1\}} < \lambda_1 \delta^2$, we have

$$-2\langle \Psi(x), x \rangle - 2\langle \Phi(x), L^{-1}x \rangle$$
$$\leq -\delta^2 \|x\|_{r+1}^{r+1} + 2|\Phi(0)|\|L^{-1}\|_{r+1}\|x\|_{r+1} + 2|\Psi(0)|\|x\|_{r+1} + \gamma \|x\|_H^2$$
$$\leq \theta_2 - \theta_1 \|x\|_{r+1}^{r+1}, \qquad x \in L^{r+1}(\mathbf{m})$$

for some $\theta_1, \theta_2 > 0$. Combining this with the Itô formula for the square of the norm, we obtain

(2.11) $$d\|X_t\|_H^2 \leq (c - \theta \|X_t\|_{r+1}^{r+1})\, dt + 2\langle Q\, dW_t, X_t \rangle_H$$

for some $c, \theta > 0$. Then

$$\mu_n(\|\cdot\|_{r+1}^{r+1}) := \frac{1}{n} \int_0^n \mathbf{E}\|X_t(0)\|_{r+1}^{r+1}\, dt$$
$$\leq \frac{c}{\theta} - \frac{1}{n}\|X_n(0)\|_H^2 \leq \frac{c}{\theta}, \qquad n \geq 1.$$

Hence, to prove the tightness of $\{\mu_n\}$, it suffices to prove that $\|\cdot\|_{r+1}$ is a compact function, that is, $K_N := \{\|\cdot\|_{r+1} \leq N\}$ is relatively compact in $H$ for any $N > 0$. Since the embedding $L^{r+1}(\mathbf{m}) \subset H$ is continuous, it follows that $\|\cdot\|_Q$ is bounded on $K_N$. Moreover, since $Q$ is Hilbert–Schmidt from $L^2(\mathbf{m})$ to $H$, $\|\cdot\|_Q$ is a compact function on $H$. Therefore, $K_N$ is relatively compact in $H$.

(b) The uniqueness and full support of $\mu$. By (1.3) with $\gamma \leq 0$ and the Itô formula, there exist $\delta, \theta > 0$ such that

$$d\|X_t(x) - X_t(y)\|_H^2 \leq -\delta^2 \|X_t(x) - X_t(y)\|_{r+1}^{r+1}\, dt$$
$$\leq -\theta \|X_t(x) - X_t(y)\|_H^{r+1}\, dt, \qquad x, y \in H.$$

Thus, $\lim_{t \to \infty} \|X_t(x) - X_t(y)\|_H = 0$, $x, y \in H$. This implies that $\mu$ is the unique invariant probability measure of $P_t$.

Next, since $\mu$ is the invariant probability measure of $P_t$, by (1.5) with $\alpha := 2$,

$$(P_t 1_A(x))^2 \int_H e^{-2c(\theta,t)\|x-y\|_H^{2(3-r+\theta)/(2+\theta)}} \mu(dy)$$

$$\leq \int_H P_t 1_A(y) \mu(dy) = \mu(A), \qquad A \in \mathcal{M}.$$

Then the transition kernel $P_t(x, dy)$ is absolutely continuous w.r.t. $\mu$ so that it has a density $p_t(x, y)$. Thus, if $\operatorname{supp} \mu \neq H$, then there exist $x_0 \in H$ and $r > 0$ such that $B(x_0, r) := \{y \in H : \|x_0 - y\|_H \leq r\}$ is a null set of $\mu$. Hence, $P_t(x_0, B(x_0, r)) = 0$. Therefore, letting $X_t(x_0)$ be the solution to (1.2) with $X_0(x_0) = x_0$, we obtain

$$\mathbf{P}(\|X_t(x_0) - x_0\|_H \leq r) = 0, \qquad t > 0.$$

Since $X_t(x_0)$ is a continuous process on $H$, this implies $\mathbf{P}(\|X_0(x_0) - x_0\|_H \leq r) = 0$ which is impossible. So, $\mu$ has full support on $H$.

(c) Concentration of $\mu$. By (2.11), for $c' := (r+1)\varepsilon_0/2$ we have

(2.12)
$$de^{\varepsilon_0 \|X_t\|_H^{r+1}} + dM_t$$
$$\leq (c - \theta\|X_t\|_{r+1}^{r+1} + 2c'\|Q\|_{\mathcal{L}_{\mathrm{HS}}}^2 \|X_t\|_H^{r+1}) c' \|X_t\|_H^{r-1} e^{\varepsilon_0 \|X_t\|_H^{r+1}} dt$$

for some local martingale $M_t$. Since $\|\cdot\|_{r+1} \geq c_0 \|\cdot\|_H$ for some constant $c_0 > 0$, when $\varepsilon_0 > 0$ is small enough there exist $c_1, \theta_1 > 0$ such that

$$de^{\varepsilon_0 \|X_t\|_H^{r+1}} \leq (c_1 - \theta_1 \|X_t\|_{r+1}^{r+1} e^{\varepsilon_0 \|X_t\|_H^{r+1}}) dt + dM_t.$$

This implies

$$\mu_n(e^{\varepsilon_0 \|\cdot\|_H^{r+1}}) \leq \frac{1}{\theta_1 n} + \frac{c_1}{\theta_1}, \qquad n \geq 1.$$

Hence, $\mu(e^{\varepsilon_0 \|\cdot\|_H^{r+1}}) < \infty$ since $\mu$ is the weak limit of a subsequence of $\mu_n$.

Finally, by (2.11) we have

$$\int_0^1 P_t \|\cdot\|_{r+1}^{r+1}(x) \, dt \leq c_2(1 + \|x\|_H^2), \qquad x \in H,$$

for some $c_2 > 0$. Thus, $\mu(\|\cdot\|_{r+1}^{r+1}) \leq c_2(1 + \mu(\|\cdot\|_H^2)) < \infty$. $\square$

PROOF OF THEOREM 1.2(2). For any $p > 1$ and any nonnegative measurable function $f$ with $\mu(f^{p/(p-1)}) \leq 1$, it follows from (1.5) with $\alpha := p/(p-1)$ that

$$(P_t f(x))^{p/(p-1)} \leq (P_t f^{p/(p-1)}(y)) \exp[pc_t \|x - y\|_H^{2(3-r+\theta)/(2+\theta)}],$$

$$x, y \in H.$$

Thus,

$$(P_t f(x))^{p/(p-1)} \int_H e^{-pc_t \|x-y\|_H^{2(3-r+\theta)/(2+\theta)}} \mu(dy) \le \mu(f^{p/(p-1)}) \le 1.$$

Therefore,

$$\langle p_t(x,\cdot), f\rangle_\mu = P_t f(x) \le \left(\int_H e^{-pc_t \|x-y\|_H^{2(3-r+\theta)/(2+\theta)}} \mu(dy)\right)^{-(p-1)/p}.$$

This implies (1.6). □

PROOF OF THEOREM 1.2(3). Let $f \in L^2(\mu)$ with $\mu(f^2) = 1$. By (1.5) with $\gamma = 0$ and constants $\xi, \delta > 0$, there exists a constant $c > 0$ depending on $r$ and $\theta$ such that

$$(P_t f)^2(x) \exp\left[-\frac{c\|x-y\|_H^{2(3-r+\theta)/(2+\theta)}}{t^{(4+\theta)/(2+\theta)}}\right] \le P_t f^2(y), \qquad x, y \in H,\ t > 0.$$

Taking integration for both sides w.r.t. $\mu(dy)$, we obtain

(2.13)
$$(P_t f)^2(x) \le \frac{1}{\mu(B(0,1))} \exp\left[\frac{c(\|x\|_H + 1)^{2(3-r+\theta)/(2+\theta)}}{t^{(4+\theta)/(2+\theta)}}\right], \qquad x \in H, t > 0,$$

where $B(0,1) := \{y \in H : \|y\|_H \le 1\}$ has positive mass of $\mu$.

If $r = 1$, then by (2.13) and Theorem 1.2(1) we have

$$\int_H (P_t f)^4(x) \mu(dx)$$
$$\le \frac{1}{\mu(B(0,1))} \int_H \exp\left[\frac{c(\|x\|_H + 1)^{2(3-r+\theta)/(2+\theta)}}{t^{(4+\theta)/(2+\theta)}}\right] \mu(dx) < \infty$$

for sufficiently big $t > 0$. Thus, $P_t$ is hyperbounded, that is, $\|P_t\|_{2\to 4} < \infty$ for some $t > 0$. Since $P_t$ has transition density w.r.t. $\mu$, according to, for example, [23] it is compact in $L^2(\mu)$ for large $t > 0$. In particular, if $\gamma \le 0$, then the process is ergodic so that its generator has a spectral gap. Thus, $\|P_t - \mu\|_2 \le ce^{-\lambda t}$ for some $c > 0$ and all $t > 0$. Therefore, by a standard argument we obtain the hypercontractivity from the hyperboundedness.

If $r > 1$, then (2.12) implies

$$de^{\varepsilon_0 \|X_t\|_H^{r+1}} \le c_2 - \theta_2 \|X_t\|_H^{2r} e^{\varepsilon_0 \|X_t\|_H^{r+1}} dt + dM_t$$

for some small $\varepsilon_0 > 0$ and some $c_2, \theta_2 > 0$. Thus, letting $h(t)$ solve the equation

(2.14)
$$h'(t) = c_2 - \theta_2 \varepsilon_0^{-2r/(1+r)} h(t) \{\log h(t)\}^{2r/(r+1)},$$
$$h(0) = e^{\varepsilon_0 \|x\|_H^{r+1}},$$

we have

(2.15) $$\mathbf{E}e^{\varepsilon_0\|X_t(x)\|_H^{r+1}} \leq h(t).$$

Since $\frac{2r}{r+1} > 1$, (2.14) and (2.15) imply

(2.16) $\mathbf{E}e^{\varepsilon_0\|X_t(x)\|_H^{r+1}} \leq \exp[c_3(1+t^{-(r+1)/(r-1)})], \qquad t > 0, x \in H,$

for some constant $c_3 > 0$. Next, by (2.13) we have

(2.17) $$\begin{aligned}\|P_t f\|_\infty &= \|P_{t/2} P_{t/2} f\|_\infty \\ &\leq c_4 \sup_{x \in H} \mathbf{E} \exp\left[\frac{c_4}{t^{(4+\theta)/(2+\theta)}}\|X_{t/2}(x)\|_H^{2(3-r+\theta)/(2+\theta)}\right],\end{aligned}$$
$$t > 0,$$

for some $c_4 > 0$. Since there exists $c_5 > 0$ such that

$$\frac{c_4}{t^{(4+\theta)/(2+\theta)}} u^{2(3-r+\theta)/(2+\theta)} \leq \varepsilon_0 u^{r+1} + c_5 t^{-(r+1)/(r-1)}, \qquad u, t > 0,$$

(1.7) follows immediately from (2.16) and (2.17). Finally, according to [23] (see also [9], Lemma 3.1), the compactness of $P_t$ follows immediately since $P_t$ is uniform integrable in $L^2(\mu)$ and has transition density w.r.t. $\mu$. $\square$

**3. Examples.** As explained in Remark 1.1, for $L := \Delta$ the Dirichlet Laplace operator, our results only apply to a space of dimension less than 2. The aim of this section is to show that, by means of spectral representation, we have many more choices of $L$ to illustrate our theorems.

Let $L_0$ be a self-adjoint operator on $L^2(\mathbf{m})$ with discrete spectrum

$$(0 \leq) \lambda_1^{(0)} \leq \lambda_2^{(0)} \leq \cdots$$

and the corresponding unit eigenfunctions $\{e_i\}_{i \geq 1}$. As in Corollary 1.3, let $Qe_i := q_i e_i$ for a sequence $\{q_i \neq 0\}_{i \geq 1}$. Let, for simplicity, $\Phi(s) = -c_0 s$ and $\Psi \in C(\mathbf{R})$ satisfy

(3.1) $$\begin{aligned}(\Psi(s_1) - \Psi(s_2))(s_1 - s_2) &\geq \delta^2 |s_1 - s_2|^{r+1}, \\ |\Psi(s)| &\leq c(1+|s|^r), \qquad s, s_1, s_2 \in \mathbf{R},\end{aligned}$$

for some $c_0 \geq 0$ and $c, \delta > 0$. For any positive and strictly increasing function $\varphi$ on $[0, \infty)$, we consider (1.2) for

$$L := -\varphi(-L_0) = -\sum_{i=1}^\infty \varphi(\lambda_i^{(0)}) \langle e_i, \cdot \rangle e_i.$$

That is, consider

(3.2) $$dX_t = -\{\varphi(-L_0)\Psi(X_t) + c_0 X_t\} dt + Q \, dW_t.$$

PROPOSITION 3.1. *Let $\inf_{i\geq 1} q_i^2 > 0$ and $\Phi = 0$, and let $\Psi$ satisfy* (3.1). *If $\varphi$ is strictly positive such that*

$$\sum_{i=1}^{\infty} \frac{q_i^2}{\varphi(\lambda_i^{(0)})} < \infty, \tag{3.3}$$

*then the Markov semigroup of the solution to* (3.2) *satisfies all assertions in Theorems* 1.1 *and* 1.2 *for any $\theta \in (r-3, r-1]$ and some $\xi > 0$.*

PROOF. Let $L := -\varphi(-L_0)$ whose eigenvalues are $-\lambda_i := -\varphi(\lambda_i^{(0)})$, $i \geq 1$. Obviously, all conditions in Corollary 1.3 are satisfied for the present situation. Thus, the proof is completed by Corollary 1.3. □

To conclude this paper, we present two examples where $L_0$ is either the Dirichlet Laplacian on a finite volume domain in $\mathbf{R}^d$ or the Ornstein–Uhlenbeck operator on $\mathbf{R}^d$, so that $L$ can be taken as high-order differential operators on $\mathbf{R}^d$ or on a domain.

EXAMPLE 3.2. Assume the situation of Proposition 3.1 but simply take $q_i = 1$, $i \geq 1$.

(1) Let $L_0 := \Delta - x \cdot \nabla$ and let $\mathbf{m}$ be the standard Gaussian measure on $E := \mathbf{R}^d$. It is well known that the set of eigenvalues of $-L_0$ is $\mathbf{Z}_+$, and the eigenspace of each $k \geq 0$ is

$$\mathrm{span}\biggl\{\prod_{i=1}^{d} H_{k_i}(x_i) : k_1 + \cdots + k_d = k, k_1, \ldots, k_d \geq 0\biggr\},$$

where $H_0 \equiv 1$ and

$$H_n(s) := \frac{(-1)^n}{\sqrt{n!}} e^{s^2/2} \frac{d^n}{ds^n} e^{-s^2/2}, \qquad s \in \mathbf{R}, \ n \geq 1.$$

Thus, there exists $\sigma > 0$ such that

$$\lambda_i^{(0)} \geq \sigma(i-1)^{1/d}, \qquad i \geq 1.$$

Then (3.3) holds for $\varphi(s) := (\varepsilon + s)^q$ for any $\varepsilon > 0$ and $q > d$, so that all assertions in Theorems 1.1 and 1.2 hold for the solution to (3.2).

(2) Let $L_0 := \Delta$ be the Dirichlet Laplace operator on a domain $D \subset \mathbf{R}^d$ with finite volume, and let $\mathbf{m}$ be the normalized volume measure on $D$. By the Sobolev inequality we have (see [21], Corollaries 1.1 and 3.1)

$$\lambda_i^{(0)} \geq \sigma i^{2/d}, \qquad i \geq 1,$$

for some $\sigma > 0$. Then (3.3) holds for $\varphi(s) := s^q$ for any $q > d/2$, so that all assertions in Theorems 1.1 and 1.2 hold for the solution to (3.2).

## APPENDIX: EXISTENCE AND UNIQUENESS OF SOLUTIONS

We first recall the following result due to [13], then derive the existence and the uniqueness for the solution to generalized stochastic porous media equations.

THEOREM A.1 ([13], Theorems II.2.1, II.2.2). *Let $H$ be a real separable Hilbert space and $V$ and $V^*$ two real Banach spaces such that the embeddings $V \subset H \subset V^*$ are dense and continuous. Let $\mathcal{L}_{\mathrm{HS}}$ be the space of all Hilbert–Schmidt operators from some real separable Hilbert space $G$ to $H$ and $W_t$ the cylindrical Brownian motion on $G$. Let $T > 0$ be fixed and*

$$A : [0,T] \times V \times \Omega \to V^* \quad \text{and} \quad Q : [0,T] \times V \times \Omega \to \mathcal{L}_{\mathrm{HS}}$$

*be progressively measurable such that*

(A1) *Semicontinuity of $A$: for any $v_1, v_2, v \in V$ and any $t \in [0,T]$, $\mathbf{R} \ni \lambda \mapsto {}_{V^*}\langle A(t, v_1 + \lambda v_2), v\rangle_V$ is continuous, where ${}_{V^*}\langle \cdot, \cdot \rangle_V$ is the duality between $V^*$ and $V$.*

(A2) *Monotonicity of $(A, Q)$: there exists a constant $K > 0$ such that for any $t \in [0,T]$,*

$$2{}_{V^*}\langle A(t, v_1) - A(t, v_2), v_1 - v_2\rangle_V + \|Q(t, v_1) - Q(t, v_2)\|^2_{\mathcal{L}_{\mathrm{HS}}}$$
$$\leq K\|v_1 - v_2\|^2_H, \qquad v_1, v_2 \in V.$$

(A3) *Coercivity of $(A, Q)$: there exist two constants $\alpha, K > 0$ and a positive adapted process $f \in L^1([0,T] \times \Omega; dt \times \mathbf{P})$ such that*

$$2{}_{V^*}\langle A(t, v), v\rangle_V + \|Q(t, v)\|^2_{\mathcal{L}_{\mathrm{HS}}} + \alpha\|v\|^{r+1}_V \leq f_t + K\|v\|^2_H$$

*holds for all $t \in [0,T]$, $v \in V$.*

(A4) *Boundedness of $A$: there exist a constant $K > 0$ and a positive adapted process $f \in L^1([0,T] \times \Omega; dt \times \mathbf{P})$ such that*

$$\|A(t,v)\|_{V^*} \leq f_t^{r/(r+1)} + K\|v\|^r_H, \qquad t \in [0,T],\ v \in V.$$

*Then for any $X_0 \in L^2(\Omega \to H; \mathcal{F}_0; \mathbf{P})$, (A.1) has a unique solution $\{X_t\}_{t \in [0,T]}$ which is an adapted continuous process on $H$ such that $\mathbf{E} \int_0^T \|X_t\|^{r+1}_V dt < \infty$ and*

$$\langle X_t, v\rangle_H = \langle X_0, v\rangle_H + \int_0^t {}_{V^*}\langle A(s, X_s), v\rangle_V\, ds + \int_0^t \langle Q(s, X_s)\, dW_s, v\rangle_H$$

*holds for all $v \in V$, $t \in [0,T]$.*

We now return to the framework in Section 1 and consider the following equation which is even more general than (1.2):

(A.1) $$dX_t = \{L\Psi(t, X_t) + \Phi(t, X_t)\}\, dt + Q(t, X_t)\, dW_t,$$

where
$$Q:[0,\infty) \times H \times \Omega \to \mathcal{L}_{HS}$$
is a progressively measurable mapping such that

(A.2)
$$\|Q(t,x)\|^2_{\mathcal{L}_{HS}} \leq h_t(1 + \|x\|^2_H),$$
$$\|Q(t,x) - Q(t,y)\|^2_{\mathcal{L}_{HS}} \leq h_t\|x - y\|^2_H$$

holds for some positive function $h \in C([0,\infty))$ and all $x, y \in H$.

THEOREM A.2. *Assume* (1.3) *and* (A.2) *for some positive function* $h \in C([0,\infty))$ *and all* $x, y \in H$.
  (1) (A.1) *has a unique solution for any* $X_0 \in L^2(\Omega \to H; \mathcal{F}_0; \mathbf{P})$.
  (2) *Let* $X_t$ *solve* (A.1) *for* $X_0 = x \in H$. *Then* (2.2) *has a unique solution*.

PROOF. (1) Let
$$A(t,x) := L\Psi(t,x) + \Phi(t,x), \qquad t \geq 0,\ x \in L^{r+1}(\mathbf{m}).$$
To make this quantity meaningful, let $V := L^{r+1}(\mathbf{m})$. Then the embedding $V \subset H$ is continuous. Let $V^*$ be the dual space of $V$ w.r.t. $H$. By (1.3) and the assumption of $L$, that is, $L^{-1}$ is bounded in $L^{r+1}(\mathbf{m})$ if $|\Phi(t,s)| \leq \sigma_t(1 + |s|^{(1+r)/2})$ does not hold for any positive $\sigma \in C([0,\infty))$, we conclude that $A(t,x)$ is well defined as an element in $V^*$ by letting
$$_{V^*}\langle A(t,x), v\rangle_V := -\langle \Psi(t,x), v\rangle - \langle \Phi(t,x), L^{-1}v\rangle, \qquad v \in V.$$

It is now easy to see that under (1.3), (A.2) and the continuity of $\Psi(t,s)$ and $\Phi(t,s)$ in $s$, all assumptions in the above theorem hold. Therefore, the proof is completed.

(2) By (1) we only have to prove (A1)–(A4) for $Q = 0$ and
$$A(t,x) := \frac{X_t - x}{\|X_t - x\|^\varepsilon_H} 1_{\{X_t \neq x\}}$$
for $\varepsilon \in (0, \frac{1}{2}]$. Since by (1.3) and the Itô formula (see [13], Theorem I.3.2) one has
$$d\|X_t\|^2_H \leq 2\langle Q(t,X_t)\,dW_t, X_t\rangle_H$$
$$- \sigma\|X_t(x)\|^{r+1}_{r+1}\,dt + (c + \|Q(t,X_t)\|^2_{\mathcal{L}_{HS}})\,dt$$
for some $c, \sigma > 0$, it follows from (A.2) that
$$\sup_{t \in [0,T]} \mathbf{E}\|X_t\|^2_H < \infty.$$

Thus, $A(t,x) \in H$ with
$$\|A(t,x)\|_H = \|X_t - x\|^{1-\varepsilon}_H, \qquad x \in H.$$

Therefore, (A1), (A3) and (A4) hold. To verify (A2), it suffices to prove

(A.3) $\quad \langle A(t,x) - A(t,y), x - y \rangle_H \leq 0 \quad \text{on } \Omega, x, y \in H.$

Without loss of generality, for a fixed $\omega \in \Omega$ we only verify (A.3) for $x, y \in H$ with

(A.4) $\quad \|X_t - x\|_H \leq \|X_t - y\|_H.$

We now prove (A.3) for the following two situations, respectively.

(i) If $\|X_t - x\|_H \geq \|x - y\|_H$, then by (A.3), the mean valued theorem and the triangle inequality, we have

$$\langle A(t,x) - A(t,y), x - y \rangle_H$$
$$= -\frac{\|x-y\|_H^2}{\|X_t - x\|_H^\varepsilon} + \frac{\|X_t - y\|_H^\varepsilon - \|X_t - x\|_H^\varepsilon}{\|X_t - y\|_H^\varepsilon \|X_t - x\|_H^\varepsilon} \langle X_t - y, x - y \rangle_H$$
$$\leq -\frac{\|x-y\|_H^2}{\|X_t - x\|_H^\varepsilon} + \frac{\varepsilon \|X_t - y\|_H^{1-\varepsilon} \|x - y\|_H^2}{\|X_t - x\|_H}$$
$$\leq -\frac{\|x-y\|_H^2}{\|X_t - x\|_H^\varepsilon} + \frac{\varepsilon(\|X_t - x\|_H^{1-\varepsilon} + \|x - y\|_H^{1-\varepsilon}) \|x - y\|_H^2}{\|X_t - x\|_H}$$
$$\leq -\frac{(1 - 2\varepsilon)\|x-y\|_H^2}{\|X_t - x\|_H^\varepsilon} \leq 0.$$

(ii) If $\|X_t - x\|_H \leq \|x - y\|_H$, then by (A.3) and the triangle inequality, we have

$$\langle A(t,x) - A(t,y), x - y \rangle_H$$
$$= -\frac{\|x-y\|_H^2}{\|X_t - y\|_H^\varepsilon} + \frac{\|X_t - x\|_H^\varepsilon - \|X_t - y\|_H^\varepsilon}{\|X_t - y\|_H^\varepsilon \|X_t - x\|_H^\varepsilon} \langle X_t - x, x - y \rangle_H$$
$$\leq -\frac{\|x-y\|_H^2}{\|X_t - x\|_H^\varepsilon} + \frac{\|x-y\|_H^\varepsilon \|X_t - x\|_H \|x - y\|_H}{\|X_t - x\|_H^\varepsilon \|X_t - y\|_H^\varepsilon}$$
$$\leq -\frac{\|x-y\|_H^2}{\|X_t - x\|_H^\varepsilon} + \frac{\|x-y\|_H^{1+\varepsilon} \|X_t - x\|_H^{1-\varepsilon}}{\|X_t - y\|_H^\varepsilon} \leq 0. \qquad \square$$

**Acknowledgments.** The author would like to thank the referees for useful comments and Mr. Wei Liu for an observation extending the earlier version to $r \geq 3$.

SCHOOL OF MATHEMATICAL SCIENCES
BEIJING NORMAL UNIVERSITY
BEIJING 100875
CHINA
E-MAIL: wangfy@bnu.edu.cn
URL: http://math.bnu.edu.cn/~wangfy/